\documentclass[12pt]{article}       
\usepackage{a4}
\usepackage{amsmath,amssymb,amsthm} 
\usepackage[all]{xy}
\theoremstyle{plain}               
\newtheorem{defn}{Definition}[section]

\newtheorem{proposition}[defn]{Proposition}
\newtheorem{corollary}[defn]{Corollary}

\theoremstyle{definition}          

\newtheorem{example}[defn]{Example}
\begin{document}
\begin{center}
{\huge{The planar Chain Rule and the Differential Equation for the planar Logarithms}
\\
\bigskip

 \Large{L. Gerritzen}}

\bigskip

 (29.09.2004)
\end{center}

\
\bigskip

\begin{abstract}
A planar monomial is by definition an isomorphism class of a finite, planar, reduced rooted tree. If $x$ denotes the tree
with a single vertex, any planar monomial is a non-associative product in $x$ relative to $m-$array grafting.
A planar power series $f(x)$ over a field $K$ in $x$ is an infinite sum of $K-$multiples of planar monomials including
the unit monomial represented by the empty tree.\\
For every planar power series $f(x)$ there is a universal differential $d f(x)$ which is a planar power series in $x$ and
a planar polynomial in a variable $y$ which is the differential $d x$ of $x$.\\
We state a planar chain rule and apply it to prove that the derivative $\frac{d}{dx}(Exp_k(x))$ is the $k-$ary planar exponential
series.
A special case of the planar chain rule is proved and it derived that the planar universe series $Log_k (1+x)$ of
$Exp_k(x)$ satisfies the differential equation
$$\biggl((1+x) \frac{d}{dx}\biggl)( Log_k(1+x)) = 1$$
where $(1+x) \frac {d}{dx}$ is the derivative which when applied to $x$ results in $1 + x.$
\end{abstract}

\bigskip

\textbf{Introduction}

\medskip

For every planar power series $f(x)$ in one variable $x$, there is a universal differential  $d f (x)$ which is a planar power
series in $x$ and a planar polynomial in a variable $y$ which is  the differential $dx$ of $x$.\\
If $g (x)$ is a planar power series of order $\geq 1$, then the substitution of $g (x)$ for $x$ in $f (x)$ gives a
power series $f (g(x))$ and the differential $d f(g(x))$ of $f(g(x))$ can be computed as
$$d \varphi_g (d f)$$
where $\varphi_g$ is the substitution homomorphism induced by $g$ and $d\varphi_g$ is the homomorphism on the algebra of
universal differential forms in $x$ and $d x$ extending $\varphi_g$ and mapping $dx$ into $dg.$\\
This formula is called the planar chain rule.

\bigskip

We note an application to the planar exponential series $E x p_k (x), k \in \mathbb{N}, k \geq 2,$ see $[G1]$, section 3.\\
One obtains a functional equation for the differential $\omega_k$ of $E x p_k (x)$, see section 6. It follows that the derivative
$\frac{d}{dx} (Exp_k(x))$ is equal to $Exp_k(x).$\\
This gives a conceptual proof of this differential relation, while the proof in [G1], section 3, is quite technical.\\
We also prove a special chain rule, if $g^\prime (x)$ where $g^\prime (x)= \frac{d}{dx}\  g(x)$ is the derivative with
respect to $x$. From it we find that
$$\biggl((1 + x) \frac{d}{dx}\biggl) ( Log_k(1+x))=1,$$
if $Log_k\ (1+x)$ is the $k-$ary planar logarithm, see [G 2], section 5.
\bigskip

In section $1$ notions about $\{x, y \}$- labeled finite, planar, reduced rooted trees are fixed. Some basic facts about the algebra
$ K \{\{x, y\}$ of power series in $x$ and polynomials in $y$ over  a field $K$ are presented in section $2$. The universal
differential
$$d: K \{\{x\}\} \rightarrow K \{\{x, y\}$$
is introduced in section $4$ and the planar chain rule is proved in section 5.\\
The application to the derivative of $Exp_k(x)$ relative to $x$ is given in section $6$.\\
In section 7 we prove a special rule for the derivative $\frac{d}{dx},$ if the substitution series $g(x)$ satisfies
the equation $\frac{d}{dx}\ g (x)=1 + g (x)$ and we show that $$\biggl(( 1 + x) \frac{d}{dx}\biggl)\  (Log_k(1 + x))=1.$$
If $h_n$ is the homogeneous component of $Log_k (1 + x)$ of degree $n$, then
$$h^\prime_{n+1}(x) =( -n) h_n(x)$$
for $n \geq 1.$\\
In [G2], section 6, the homogeneous component $h_4$ of $Log_k(1+x)$ has been computed to be
$$h_4(x)= \biggl(\frac{k-3}{3![2]} -\frac{1}{4![3]!} (k+1) (k-2) (k-3)\biggl) \cdot x^4 +$$
\medskip
$$+\biggl(\frac{1}{2} \frac{1}{3![2]}- \frac{2}{4![3]!}(k-2)\biggl) (x \cdot x^3 + x^3 \cdot x)$$
\medskip
$+\biggl(\frac{1}{3![2]!} \cdot \frac{3}{2} - \frac{1} {8} - \frac{3}{4![3]!}\biggl) (x\cdot (x \cdot x^2) +
 x\cdot(x^2 \cdot x) + (x\cdot x^2) \cdot x +(x^2 \cdot x) \cdot x)$
\medskip
$$+\biggl(\frac{1}{3![2]!} \cdot \frac{3}{2} - \frac{1}{8} - \frac{3(k+1)}{4![3]!}\biggl) \cdot x^2 \cdot x^2$$
\medskip
$+\biggl(\frac{1}{2} \frac{k-2}{3![2]!}- \frac{2(k+1)(k-2)}{4![3]!}\biggl) (x \cdot x \cdot x^2+ x \cdot x^2 \cdot x +
x^2 \cdot x \cdot x).$

\bigskip

From this it follows that
$$ h_3(x)=(-\frac{1}{3}) h^\prime_4(x)= \frac{1}{4} \frac{k}{[2]}(x \cdot x^2 + x^2 \cdot x) -\frac{1}{3!} \frac{k-2}{[2]}\cdot x^3,$$
$$ h_2(x)=(-\frac{1}{2}) h^{\prime}_3(x) = -\frac{1}{2}\ x^2,$$
$$ h_1(x) = - h^\prime_2(x)= x.$$
\section{Labeled planar rooted trees}

\bigskip

Let $x, y$ be different elements; they will be called variables in the sequel.\\
$A\  \{x,y\}$- labeled finite, planar reduced rooted tree is a pair $S=(T,\lambda)$ where $T$ is a finite, planar, reduced
rooted tree, see [G1], and $\lambda$ is a map $L(T) \rightarrow \{x, y\}$ where $L (T)$ denotes the set of leaves of $T$. \\
Then
$$deg_x(T, \lambda):= \# \lambda^{-1}(x)$$
is called the $x$-degree of $(T, \lambda)$ which is the number of leaves of $T$ labeled with $x$ by $\lambda$.\\
Also
$$ deg_y(T, \lambda):= \# \lambda^{-1}(y)$$
is called the $y$-degree of $(T, \lambda)$
which is the number of leaves of $T$ labeled with $y$ by $\lambda$.
Then\ $ deg\  (T, \lambda)=deg\  T:= \# L(T)$ is called the total degree of $(T, \lambda)$.

\bigskip

We denote by $P(x,y)$ the set of isomorphism classes of $\{x,y\}-$ labeled finite, planar, reduced rooted trees.\\
Given $S_1,..., S_m \in P(x,y), S_i = (T_i, \lambda_i)$ for $ 1 \leq i \leq m$. Then there is a unique $S \in P (x,y)$
with the following properties:\\ $S=(T, \lambda)$ and if $ \rho_T$ is the root of $T$, then $T - \rho_T$ is a planar rooted\\
forest whose $i$- th component is isomorphic with $T_i$. The set of leaves of $L(T)$ is the disjoint union of $L(T_1),...,
L(T_n)$ and the labeling $\lambda$ induces $\lambda_i$ on $T_i$ for all $ 1 \leq i \leq m.$\\
We recall that a planar rooted forest is an ordered disjoint system of planar rooted trees.\\
The construction described above gives a unique map
$$ \bullet_m: (P(x,y))^m \rightarrow P(x,y)$$
which is called the $m$-ary grafting operation for labeled planar rooted trees.\\
Let $P^\prime (x, y) = P(x,y) \cup \{1_P\}$ where $1_p$ is a symbol representing the empty tree. There is a unique extension
of the $m-$ ary grafting operation to a map
$$ \bullet_m: (P^\prime (x,y))^m \rightarrow P^\prime (x,y)$$ satisfying the properties: $\cdot_m (1_P,...,1_P)= 1_P$\\
$\cdot_2(S, 1_P) = \cdot_2 (1_P, S)= S$ for all $S\in P^\prime (x,y),\\
\cdot_m (S_1,...,S_m)= \cdot_{m-1} (S_1,...,
S_{i-1}, S_{i+1},...,S_m)$ if $S_i = 1_P$ for all $S_i \in P^\prime(x,y).$\\
It is also refered to as the $m$-ary grafting operation.

\section{Polynomials and Power series}

Let $K$ be a field and $V$ be the $K$- vectorspace of all $K-$ valued functions $f$ on  $P^\prime
(x,y)$ such that for all $n \in \mathbb{N}$ the set
$$\{S \in P^\prime (x,y):f(S)\not= 0, ord_x(S)=n\}$$ is finite.
For any $f \in V$ the value $f(S)$ at $S \in P^\prime(x,y)$ will be denoted by $c_S(f).$ It will be called the coefficient of
$f$ relative to $S$.\\
For any $m \in \mathbb{N}, m \geq 2,$ we are going to define a $m-$ary operation $\bullet_m$ on $V$ by putting
$$ f:= \bullet_m(f_1,...,f_m)$$
if  $c_S(f) :=  \sum c _{S_1}(f_1)\cdot...\cdot c_{S_m}(f_m)$ where the summation is extended over all $S_1, ..., S_m$ such that
$S= \bullet_m (S_1, ... , S_m).$

It is easy to check that $\bullet_m$ is a $K-$ multilinear map $V^m \rightarrow V.$
We will write\\ $f_1 \times f_2 \times...\times f_m$ for $\bullet_m(f_1,...,f_m)$ sometimes.

\bigskip

In this note a $K-$algebra is by definition a $K-$ vectorspace $W$ together with a family $(\bullet_m)_{m \geq 2}$ of
$K-$ multilinear maps $\bullet_m:W^m \rightarrow W.$ In the context of operad theory such an algebra can be defined as
an algebra over an operad in the category of $K-$ linear spaces which is freely generated by $\{\mu_m \colon  m \geq 2\}$
where $\mu_m$ has degree $m$, see [MSS], Chap. (1.9). It is a tree operad, see [BV].
Trees are inspired by the attempt to obtain a general composite operation from a collection of indecomposable operations,
see [BV].

\bigskip

Other authors, see [SV], use the term hyperalgebra to denote this type of algebras.
For the general theory of operads and algebras over operads see also [L], [LR], [H], [F].

Therefore $(V,(\bullet_m)_{m\geq 2)}$ is a $K-$algebra; it is denoted by $K\{\{x,y\}$ and will be refered to as the algebra of
planar power series in $x$ and polynomials in $y$.\\
For any $S \in P^\prime (x,y)$ there is a unique $f_S \in K\{\{x,y\}$ such that
\begin{equation*}c_T(f_s)=
\begin{cases}
1: T=S\\
0: T \in P^\prime(x,y), T \not= S.
\end{cases}
\end{equation*}
We identify $f_S$ with $S$ and thus consider $P^\prime(x,y)$ as a subset of $K\{\{x,y\}.$ It is easy to check that
$P^\prime(x,y)$ is a $K-$linearly independent subset of $K\{\{x,y\}$ and that the $K-$ vectorspace of $K\{\{x,y\}$
generated by $P^\prime(x,y)$ is a subalgebra. It is denoted by $K\{\{x,y\}$ and is called the algebra of planar polynomials:
in $x$ and $y$ over $K$.\\
For $f \in K\{\{x,y\}$ define $$ord_x(f)=\infty,\textrm {if}\  f=0$$ and if $$f \not= 0,$$ then $$ord_x(f)=min\{
n \in \mathbb{N}:\textrm {there is } S \in P(x,y)\ \textrm{with}\  deg_x(S)=\ n\  \textrm{and}\  c_S(f) \not= 0.\}$$
It is called the $x-$order of $f$\\
Let
\begin{equation*}|f|_x:=
\begin{cases}
0: f = 0\\
(\frac{1}{2})^{ord_x}(f):f \not= 0.
\end{cases}
\end{equation*}
Then $|f_1 \cdot f_2\cdot...\cdot f_m|_x= |f_1|_x \cdot|f_2|_x\cdot...\cdot|f_m|_x$\\
as $ord(f_1 \cdot f_2 \cdot...\cdot f_m)= \displaystyle\sum_{i=1}^m ord(f_i)$ and $|\  ,\  |_x$ defines a distance function on
$K\{\{x,y\}$ by $|f,g|_x: = |f - g|_x$ which is called the $x-$ adic distance.

\begin{proposition}
Any Cauchysequence in $K\{\{x,y\}$ relative to the $x-$ adic
distance is convergent.
\end{proposition}

\begin{proof}
By standard arguments.
\end{proof}

\section{Substitution homomorphisms}

\begin{proposition}

Let $g, h \in K\{\{x,y\}$ and let $ord_x(g) \geq 1.$\\
Then there is a unique $K-$algebra homomorphism
$$ \varphi_{(g,h)} = \varphi: K\{\{x, y\} \rightarrow K\{\{x,y\}$$
such that
\begin{itemize}
\item[(i)]
$\varphi(x) = g$
\item [(ii)]
$\varphi(y) = h$
\item [(iii)]
$\varphi$ is continuous with respect to the $x$-adic topology.
\end{itemize}
It is called substitution homomorphism induced by ( g, h).\\
\end{proposition}

\begin{proof}1)
First we define $\varphi^\prime (S)$ for any $S \in P(x,y).$ If $S=(T, \lambda)$ where $T$ is a planar rooted tree of
degree , there is a $m-$ary operation $\cdot_T$ associated to $T,$ see [G2], Proposition (2.7)\\
Let $\varphi^\prime(S):= \cdot_T(Z_1,...,Z_m)$ where $Z_i = g$ if $\lambda (l_i)= x$ and $Z_i= h$ if $\lambda(l_i)
= y.$ Here $l_i$ is the $i$-th leaf of $T$.\\
2 )There is a unique $K$-linear extension of the map $\varphi^\prime$ in 1) to a linear map
$$ \varphi: K\{\{x, y\} \rightarrow K\{\{x, y\}.$$
3) If $f \in K\{\{x, y\}, f_n:= \sum c_S(f) \cdot S$ where the sum is extended over all $S \in P(x, y)$ for which $
deg_x(S)= n,$ then $f_n$ is a finite sum and
$$ f= \displaystyle \sum^\infty_{n=0}f_n.$$
 From 1) it follows that $ord_x(\varphi(f)) \geq ord_x (f)$ for all $f \in K \{\{x, y\}$ as $ord_x(g) \geq 1.$ Thus
$ord_x (\varphi (f_n)) \geq n$ for all $n$.\\
One defines $\varphi (f) = \displaystyle \sum^\infty_{n=0} \varphi (f_n).$ The infinite sum is converging, as the sequence
of partial sums is a Cauchysequence with respect to $|\  |_x$.\\
Thus a continuous $K$- linear map $\varphi: K\{\{x, y\}\  K\{\{x, y\}$ is well defined. One can check that $\varphi$ is an algebra
homomorphism because its restriction to $K\{\{x, y\}$ is an algebra homomorphism by the definition in 1) and 2).\\
\qed

\bigskip

\textbf{Notation:} For $\varphi = \varphi_(g,h)$ and $f \in K\{\{x, y\}$, we denote $\varphi(f)$ also by $f(g,h).$ Then
$f = f(x,y)$ as $\varphi(x,y)=i d.$ Also $f (g, h)= f(g(x,y), h(x,y)).$ It is obtained by replacing $x$ by $g$ and $y$
by $h$.\\
For any $k \in K$ we get $f(k x, k y)= k^m f(x, y)$ , if $f$ is homogeneous of degree $m$.\\
\section{ Universal derivation}
We denote by $K\{\{x\}\}$ the closed unital $K$-subalgebra of $K\{\{x, y\}$ generated by $x$.\\
\begin{proposition}
There is a unique continuous $K$-linear map
$$d:K \{\{x\}\} \rightarrow K\{\{x, y\}$$
such that
\begin{itemize}
\item[(i)]
$d(x) = y$
\item[(ii)]
$d (f_1 \cdot f_2 \cdot ...\cdot f_m)=\\
=\displaystyle \sum^{m}_{i=1} (f_1 \cdot ... \cdot f_{i-1} \cdot  d(f_i) \cdot f_{i+1} \cdot ... \cdot  f_m)$ for all
$m \geq 2$ and all $f_1, ...,  f_m \in K\{\{x, y\}$
\end{itemize}
The map $d$ is called the universal derivation on $K\{\{x\}\}.$ It is customary to denote $d(f)$ by $d f$ and call
it the differential of $f$. Especially $y= dx.$
\end{proposition}
\begin{proof}1) First $d^\prime(S)$ is defined for all $S \in P(x, y).$
If $S = 1_P,$ then $d^\prime(1_P):=0$. If $S \not= 1_P, S=(T, \lambda), \lambda (l_i) = x$
for all leaves  $l_1, ... ,l_m$  of $T, m = deg (T), L (T) = \{l_1, ... , l_m\}, then$\\
let $d^\prime(S) := \displaystyle \sum^m_{i=1} (T, \lambda^{(i)}),$ where $\lambda^{(i)}(l_j)= y$ if $j=i$ and
$\lambda^{(i)}(l_j) = x$ if $j \not= i$.\\
2) As $P^\prime(x)$ is $K$-linearly independent, the map $d^\prime $ in 1) has a unique extension to a $K-$linear map
$$ d_0: K\{x\} \rightarrow K\{x, y\}.$$
3) There is a unique continuous extension of $d_0$ to a $K-$ linear map
$$ d: K\{\{x\}\} \rightarrow K \{\{x, y\}.$$
It is easy to check that property (ii) holds.
\end{proof}
\section{ Chain rule}
Let $g \in K \{\{x\}\}, ord(g) \geq 1.$ As a special case of Proposition (3.1) we get a continuous $K-$algebra homomorphism
$$ \varphi_g=\varphi: K \{\{x\}\} \rightarrow K \{\{x\}\}$$
such that $\varphi(x)=g.$\\
Also there is a unique continuous $K-$algebra homomorphism
$$ dy: K\{\{ x, y\} \rightarrow K \{\{x, y\}$$
such that $$d \varphi(x) = \varphi(x) = g$$
$$d\varphi(y) = dg$$
\begin{proposition}
$$ (d \varphi) \circ d = d \circ \varphi.$$
\end{proposition}
\begin{proof}
1) Let $\delta:= (d \varphi) \circ d - d \circ \varphi.$ Then $\delta$ is a continuous $K-$ linear map from
$K\{\{x\}\}$ into $K\{\{x , y\}\}$ which satisfies
$$\delta(x) =0$$
2) It is easy to see that
$$ \delta(f_1 \cdot ... \cdot f_m)= \displaystyle \sum^{m}_{i=1} \varphi(f_1)\cdot ...\cdot \varphi(f_{i-1}) \cdot
\delta(fi)\cdot \varphi (f_{i+1}) \cdot ... \cdot \varphi(f_m.$$
Using this formula one can prove by induction on the degree that for any $S \in P(x)$ one gets $\delta (S) = 0.$\\
If $S = S_1 \cdot ... \cdot S_m,$ then $deg_x(S_i) < deg_x(S)$ and by induction hypothesis $\delta(S_i)= 0, (\forall i?)$ By the
formula above then $\delta(S) = 0$\\
3) As $\delta$ is continuous and $K\{x\}$ is dense in $K\{\{x\}\}$ and $\delta(f) = 0$ by 2) and the linearity of
$\delta$, we see that $\delta= 0$.
\end{proof}

\bigskip

\textbf{Example:}\ Let $$g = k  \cdot x, k \in K.$$ Then
$$\varphi(f(x)) = f(kx)$$ and $$(d\varphi) \omega (x,y)= \omega(k x, ky).$$
If $f$ (resp. $\omega$) is homogeneous of degree $m$, then $f(kx) = k^mf(x), \omega(kx, ky) = k^m \omega (x, y).$
\end{proof}

\section{The derivative of the exponentials}
Let now $char(K)= 0, k \in \mathbb{N}k \geq 2,$ and $f= Exp_k(x)$ be the $k-$ary exponential series, see [G1].\\
Let $\omega = df$ be the differential of $f$. Then $\omega(k x, k y)= d f^k$ by the chain rule.\\
If $f= \displaystyle \sum^\infty_{n=\circ} f_n$ and $f_n$ is the homogeneous part of degree $n$, then
$$ \omega = \displaystyle \sum^\infty_{n=1} df_n$$
and $df_n$ is the homogeneous component $\omega_n$ of $\omega$ of degree $n$ with $df_0 = 0.$\\
More precisely:
$$ deg_x(df_n) = n-1,$$
$$ deg_y(df_n)= 1.$$
Now by the functional equation satisfied by $f$, we get
$$ f^k(x) = f(kx) = \displaystyle \sum^\infty_{n=0} k^nf_n(x)$$
$(\ast)$    $$ k^n f_n(x) = \displaystyle \sum_{i_1+...+i_k=n} f_{i_1} \cdot f_{i_2} \cdot ... \cdot f_{i_k}$$
where the summation is over all
$$(i_1,..., i_k) \in \mathbb{N}^k$$
with
$$ i_1 + ... + i_k=n.$$
\begin{proposition}
$$k^n \omega_n(x,y) = \displaystyle \sum^k_{j=1} \displaystyle\sum_{i_1 + ... +i_k=n} fi_1 \cdot ... \cdot fi_{j-1}
\cdot \omega (fi_j)\cdot  fi_{j+1} \cdot ... \cdot fi_k.$$
\end{proposition}
\begin{proof}
Immediate by property (ii) of Proposition 4.1.
\end{proof}

\bigskip

There is a unique continuous unital algebra homomorphism
$$\psi: K\{\{x, y\} \rightarrow K \{\{x\}\}$$
such that $\psi(x) = x, \psi(y)= 1.$\\
Recall that $\psi(f) = f$ for all $f \in K \{\{x\}\}.$\\
Then $(f \circ d) : K \{\{x\}\} \rightarrow K\{\{x\}\}$ is the \\
derivative relative to $x$, also denoted symbolically by $\frac{d}{dx},$ see [G1]
\begin{proposition}
$\frac{d}{dx} (Exp_k(x))= Exp_k(x).$
\end{proposition}
\begin{proof}1)
We have to show that the derivative $f_n^\prime = \frac{d}{dx}(f_n)$ of the homogeneous component of $f(x) = Exp_k(x)$
is equal to $f_{n-1}$ for $n \geq 1$.\\
By applying $\psi$ to the formula in Proposition (6.1) we get
$$k^nf_n^\prime(x) = \displaystyle \sum^k_{j=1}\  \displaystyle \sum_{i_1 +...+i_k=n} f_{i_1} \cdot...\cdot f_{i_j-1}
\cdot f^\prime_{i_j} \cdot f_{i_j+1} \cdot ... \cdot f_{i_k}.$$
2) We will show that $f_n^\prime(x) = f_{n-1}(x)$
for $n \geq 1$ by induction on $n$.\\
If $n= 1,$ then it is obviously true as $ord_x(f-(1+x)) \geq 2.$\\
Let now $n \geq 2.$ As $f^\prime i_j= 0$ if $i_j=0$ and $f^\prime_{i_j}= f_{i_j-1}$ for $n >0_j
\geq 1,$ we obtain\\
$$k^nf_n^\prime(x) =  \sum^k_{j=1}  \sum_{\circ 1+...yi_k=n\atop
i_j \geq 1}f_{i_1} \cdot ..., \cdot f_{i_j-1} \cdot  f_{i_j-1} \cdot f_{i_j+1} \cdot ...\cdot f_{i_k}$$
Now for any $1 \leq i \leq k$ we get
$\displaystyle \sum_{i_1 t ... tI_k=n\atop i_j \geq 1} f_{i_1} \cdot ... \cdot f_{i_j-1} \cdot f_{i_j}-1 \cdot f_{i_j+1}
\cdot ... \cdot fi_k = \displaystyle \sum_{i_1 +...+i_k = n -1} f_{i_1} \cdot  f_{i_2} \cdot ... \cdot  f_{i_k}= k^{n-1} -f_{n -1} (x)$

\bigskip

by \ \ \ \  $(\ast)$\\
It follows that
$$k^nf_n^\prime (x) = k \cdot k^{n-1} \cdot f_{ n -1}(x)$$ and
$$f^\prime_n(x) = f_{n-1}(x).$$
\end{proof}

\section{Special chain rules and logarithm}

Let $h \in K\{\{x\}\}.$ Then there is a unique algebra homomorphism $\varphi_{(x,h)}= \varphi: K\{\{x,y\} \rightarrow
K\{\{x\}\}$ such that $\varphi (x) = x$ and $\varphi(y) = h$ Then
$$\theta:= \varphi\  o\  d: K \{\{x\}\} \rightarrow
K\{\{x\}\}$$ is a derivation on $K \{\{x\}\}$ which means that is $K-$linear, continuous and satisfies the general product
rule
$$\theta(f_1 \cdot ... \cdot f_m)= \sum^m_{i=1} f_1 \cdot ... \cdot f_{i-1} \cdot \theta(f_i) \cdot f_{i+1} \cdot ... \cdot
f_m$$
for all $$f_1,..., f_m \in K\{\{x\}\}$$
This derivation $\theta$ will be formally denoted by $h \frac{d}{dx}.$\\
Be aware that $(h \frac{d}{dx}) (f) \not= h \cdot \frac{d}{dx}(f)$ in general.\\
We denote $\frac{d}{dx}(f)$ also by $f^\prime$ and call it the derivative of $f$ with respect to $x$.\\
\begin{proposition}:
Let $g \in K \{\{x\}\}$ such that $g^\prime = 1 + g$\\
Then
$$ \frac{d}{dx} (f( g(x))= f^\prime (g(x)) + (x \frac{d} {dx})(f) (g(x)) = ((1 +x) \frac{d}{dx} (f) (g(x))$$
\begin{proof}:1) First we will prove the formula in case $f$ is a tree $T$ of degree $n$\\
If $n=1,$ then $T = x$ and $\frac{d}{dx} x = 1$ while $(1+x) \frac{d}{dx}(T)= 1 + x.$ Thus $f(g(x)) = g(x)$ and $\frac
{d}{dx}(g) = 1 + g$ from which the formula follows. Let now $n > 1.$ Then
$$\frac{d}{dx}(T) = \sum^n_{i=1} T^{(i)}$$
where $T^{(i)} = \bullet_T (z_1, \cdot, z_n)$ with
$$z_j=\begin{cases}
x \colon j \not= i\\
1_P \colon j = i
\end{cases}$$
and
$$x\frac{d}{dx} (T)= n \times T$$
Also
$$\frac{d}{dx}(T(g(x)) = \sum^n_{i=1} S^i$$
with $$S^{(i)}= \bullet_T(W_1, ... , W_n)$$
and
$$w_j= \begin{cases}
g(x) \colon j \not= i\\
g^\prime(x) \colon j=i
\end{cases}$$
As $g^\prime = 1 + g,$ we obtain from the multilinearity of $\bullet_T,$ that
$$S{(i)}= \bullet_T(g, ..., g) + \bullet_{T^{(i)}} (g,... ,g)$$
From these computations it follows that the proposition holds, if $f$\ is a tree\\
2) From 1) one can extend the result to polynomials by multilinearity. As both sides of the formula are continuous
in $f$ the result follows for power series.

\end{proof}
\end{proposition}
Let $Log_k(1 + x)$ be the $k-$ary planar logarithm. see [G2], section 5.\\
\begin{proposition}:
$$((1+x) \frac{d}{dx} (Log_k( 1 + x)= 1$$
for all
$$k \in \mathbb{N},\  k > 2.$$
\end{proposition}
\begin{proof}
Let $$g= Exp_k(x) -1$$
Then
$$Log_k(Exp_k(x)) = x$$
From the special chain rule above, we can conclude that
$$(((1+ x) \frac{d}{dx})(Log_k)) Exp_k(x))=1$$
From this the formula follows as $f(g(x))=1$ only if $f(x) = 1.$
\end{proof}

\begin{corollary}:
Let $h_n$ be the homogeneous component of $Log_k (1+x)$ of degree $n.$\\
Then $$h_0 = 0$$ and
$$h^\prime_{n+1} = -n h_n$$ for
$$n \geq 1$$
\end{corollary}
\begin{proof}:
$$Log_k(1 + x)= \displaystyle \sum^\infty_{n=0} h_n =:h$$ and
$$(1 + x) \frac{d}{dx} (h_n) = h^\prime_n + n h_n,$$ as $$(s \frac{d}{dx} (h_n)= n \times h_n.$$
 Moreover $h^\prime_n$
is homogeneous of degree $(n-1).$\\
Thus $$h^\prime_{n+1} + n h_n = 0$$
for all
$$n\geq 1$$
because
 $$(1 + x) \frac{d}{dx} (h) = 1.$$
\end{proof}

\begin{example}
We give the homogeneous terms $h_n(x)$ of $Log_k(1+x)$ of degree $n \leq 4$.\\
Recall that $$[n]= \frac{k^n-1}{k-1}= \sum_{i=0}^{n-1} k^i $$
and that $$[n]!= \prod_{i=1}^n [i].$$
So especially $$[2]=k+1,\ [3] = k^2+k+1,[3]!=(k+1)(k^2+k+1)=k^3 + 2k^2+2k+1.$$

For any planar rooted tree $T$ in $P,$ let $\{T\}$ be the sum of all planar trees $S$ in $P$ whose underlying rooted tree is
isomorphic to the underlying rooted tree of $T.$\\
One can show that $Exp_k(x)$ and $Log_k(x)$ are infinite sums over $\{T\}$ with $ T \in P.$ \\
Then the homogeneous components $h_n(x)$ for $n\leq 4$ can be written as, see [G2]:

\bigskip


$h_1(x) =  x$\\
\medskip
$h_2(x) =  -\frac{1}{2}\ x^2$\\
\medskip
$h_3(x) =  \frac{1}{4} \frac{k}{[2]}\{x \cdot x^2\} -\frac{1}{3!} \frac{k-2}{[2]}\cdot x^3$\\

\bigskip

$h_4(x) = (\frac{k-3}{3![2]} -\frac{1}{4![3]!} (k+1) (k-2) (k-3)) \cdot \{x^4\} +$\\
\medskip
$+(\frac{1}{2} \frac{1}{3![2]}- \frac{2}{4![3]!}(k-2)) \{x \cdot x^3\} $\\
\medskip
$+(\frac{1}{3![2]!} \cdot \frac{3}{2} - \frac{1} {8} - \frac{3}{4![3]!})\{x \cdot(x \cdot x^2)\}$\\
\medskip
$+(\frac{1}{3![2]!} \cdot \frac{3}{2} - \frac{1}{8} - \frac{3(k+1)}{4![3]!}) \cdot \{ x^2 \cdot x^2\}$\\
\medskip
$+(\frac{1}{2} \frac{k-2}{3![2]!}- \frac{2(k+1)(k-2)}{4![3]!})+ \{x \cdot x \cdot x^2\}$
\medskip


Easy computations give

$$\frac{d}{dx}(\{x \cdot x \cdot x^2\}) = 3\  \{x \cdot x^2\} + 6\  \{x \cdot x \cdot x\}$$
$$ \frac{d}{dx} ( \{x^2 \cdot x^2 \}) = 2\  \{x \cdot x^2\}$$
$$\frac{d}{dx} ( \{x \cdot (x \cdot x^2)\} ) = 6\  \{ x \cdot x^2\} $$
$$\frac{d}{dx} \{ x \cdot x^3\} = 2 \{ x^3\} + 3\  \{x \cdot x^2\}$$

from which one can derive that
$$ h^\prime_4 (x) = - 3 h_3 (x)$$

Similarly one gets

$$h^\prime_3 (x) = - 2 h_2 (x)$$
$$ h^\prime_2 (x) = - h_1 (x)$$

In [DG] $Log_2(x)$ was computed to be $$Log_2(x)\ =\ x\ -\frac{1}{2}\ x^2\ +\ \frac{1}{3}\  (\frac{1}{2}\  x\cdot x^2\ +\ \frac{1}{2}\  x^2\cdot x)\ -$$

$- \frac{1}{4} (\frac{4}{21}\ x\cdot (x \cdot x^2)\ +\ \frac{4}{21}\  x\cdot(x^2\cdot x)\  +\frac{5}{21}\  x^2 \cdot x^2\ +\  \frac{4}{21}
(x \cdot x^2) \cdot x + \frac{4}{21} (x^2\ \cdot x)\ \cdot x)+$\\
higher terms.\\
It is obtained by substituting $2$ for $k$ in $Log_k(x)$.
\end{example}

\textbf{Open question:} In there a procedure to construct $h_{n+1}(x)$ which is an integral of $(-n) \cdot h_n(x)$ directly from
$h_n(x)?$\\
The problem arises from the fact that the space of homogeneous polynoms of degree $n+1$ whose derivative is zero is non-
trivial for $n \geq 2$

\end{document}